\documentclass[a4paper,12pt]{article} 
\usepackage[english]{babel}
\usepackage{graphicx,latexsym,amsmath,amsthm,amssymb,color,indentfirst,xfrac}

\theoremstyle{plain}
\newtheorem{thm}{Theorem}[section]
\newtheorem{cor}[thm]{Corollary}
\newtheorem{lem}[thm]{Lemma}

\theoremstyle{remark}

\theoremstyle{definition}

\newtheorem{ex}[thm]{Example}

\newcommand{\be}{\begin{equation}}           
\newcommand{\ee}{\end{equation}}

\newcommand{\ba}{\begin{align}}                
\newcommand{\ea}{\end{align}}

\newcommand{\bal}{\begin{align*}}              
\newcommand{\eal}{\end{align*}}

\newcommand{\bxx}{\begin{ex}}

\newcommand{\exx}{\end{ex}}

\newcommand{\txsit}

\newenvironment{pr}
{\begin{trivlist}
\item[\hskip\labelsep{\bf Proof.}]}                     
{$\hfill\Box$\end{trivlist}}

\title{Complexity of equitable tree-coloring problems}

\author {Keaitsuda Maneeruk Nakprasit\\ {\small\em Department of Mathematics, Faculty of Science, Khon Kaen University, 40002, Thailand }\\
{\small\em E-mail address: kmaneeruk@hotmail.com}
\and Kittikorn Nakprasit \footnote{Corresponding Author} \\
{\small\em Department of Mathematics, Faculty of Science, Khon Kaen University, 40002, Thailand }\\
{\small\em E-mail address: kitnak@hotmail.com}}
\date{}

\begin{document} 

\maketitle
\begin{abstract} 
A $(q,t)$\emph{-tree-coloring} of a graph $G$ is 
a $q$-coloring of vertices of $G$   
such that the subgraph induced by each color 
class is a forest of maximum degree at most $t.$ 
A $(q,\infty)$\emph{-tree-coloring} of a graph $G$ 
is a $q$-coloring of vertices of $G$   
such that the subgraph induced by each color 
class is a forest. 

Wu, Zhang, and Li introduced the concept of 
\emph{equitable $(q, t)$-tree-coloring} 
(respectively, \emph{equitable $(q, \infty)$-tree-coloring})  
which is a $(q,t)$-tree-coloring  
(respectively, $(q, \infty)$-tree-coloring)   
such that the sizes of any two color classes differ by at most one. 
Among other results, they obtained a sharp upper bound on 
the minimum $p$ such that $K_{n,n}$ has an equitable 
$(q, 1)$-tree-coloring for every $q\geq p.$ 

In this paper, we obtain a polynomial time criterion 
to decide if a complete bipartite graph has 
an equitable $(q,t)$-tree-coloring or 
an equitable $(q,\infty)$-tree-coloring. 
Nevertheless, deciding if a graph $G$ in general 
has  an equitable $(q,t)$-tree-coloring or 
an equitable $(q,\infty)$-tree-coloring is NP-complete. 
\end{abstract}

\section{Introduction} 

 Throughout this paper, all graphs are finite, undirected, and 
simple. We use  $V(G)$ and $E(G),$ respectively, 
to denote the vertex set and edge set of a graph $G.$ 
Let $K_{m,n}$ be a complete bipartite graph in which  
partite set $X$ has size $m$ and 
partite set $Y$ has size $n.$

An \emph{equitable $k$-coloring} of a graph is a proper vertex
$k$-coloring such that the sizes of every two color classes differ
by at most $1.$

Hajnal and Szemer\'edi~\cite{HS} settled a conjecture of Erd\H os by
proving that  every graph $G$ with maximum degree at most $\Delta$
has an equitable $k$-coloring for every $k\geq 1+\Delta.$ 
This result is now known as Hajnal-Szemer\'edi Theorem. 
Later, Kierstead and
Kostochka~\cite{KK08} gave a simpler proof of Hajnal-Szemer\'edi
Theorem. The bound of the
Hajnal-Sz{e}mer\' edi theorem is sharp, but it can be improved for
some important classes of graphs. In fact, Chen, Lih, and
Wu~\cite{CLW94} put forth the following conjecture: 
Every connected graph $G$ with maximum degree $\Delta\geq 2$ has an
equitable coloring with $\Delta$ colors, except when $G$ is a
complete graph or an odd cycle or $\Delta$ is odd and
$G=K_{\Delta,\Delta}.$

Lih and Wu~\cite{LW} proved the conjecture for bipartite graphs.
Meyer \cite{M} proved that every forest with maximum degree $\Delta$
has an equitable $k$-coloring for each $k \geq 1+\lceil
\Delta/2\rceil $ colors. This result implies the conjecture holds
for forests.  Yap and Zhang~\cite{YZ1} proved that the
conjecture holds for outerplanar graphs. Later Kostochka~\cite{Ko}
improved the result  by proving that every
outerplanar graph with maximum degree $\Delta$ has an equitable
$k$-coloring for each $k \geq 1+\lceil \Delta/2\rceil.$ 

In~\cite{ZY98}, Zhang and Yap essentially proved the conjecture
holds for planar graphs with maximum degree at least $13.$ Later
Nakprasit~\cite{Nak12} extended the result to all planar graphs with
maximum degree at least  $9.$ 
Some related results are about planar graphs without some restricted cycles 
~\cite{LiBu09, NN12, ZhuBu08}. 

Moreover, the conjecture has been confirmed for other classes of graphs, 
such as graphs with degree at most 3~\cite{CLW94, CY12}  
and series-parallel graphs \cite{ZW11}. 

A $(q,t)$\emph{-tree-coloring} of a graph $G$ is 
a $q$-coloring of vertices of $G$   
such that the subgraph induced by each color 
class is a forest of maximum degree at most $t.$ 
A $(q,\infty)$\emph{-tree-coloring} of a graph $G$ 
is a $q$-coloring of vertices of $G$   
such that the subgraph induced by each color 
class is a forest. 

In \cite{Fan11}, Fan, Kierstead, Liu, Molla, Wu, and Zhang 
considered an equitable relaxed coloring. 
They proved that every graph with maximum degree $\Delta$ has an 
$\Delta$-coloring such that each color class induces a graph 
with  maximum degree at most one and the sized of 
any two color classes differ by at least one. 

On the basis of the aforementioned research,  
Wu, Zhang, and Li \cite{WZL13}  introduced the concept of 
\emph{equitable $(q, t)$-tree-coloring} 
(respectively, \emph{equitable $(q, \infty)$-tree-coloring})  
which is a $(q,t)$-tree-coloring  
(respectively, $(q, \infty)$-tree-coloring)   
such that the sizes of any two color classes differ by at most one. 
Thus, the result of Fan, Kierstead, Liu, Molla, Wu, and Zhang can be restated that  
every graph with maximum degree $\Delta$ 
has an equitable $(\Delta,1)$-tree-coloring.  

Among other results, Wu, Zhang, and Li \cite{WZL13} 
obtained a sharp upper bound on 
the minimum $p$ such that $K_{n,n}$ has an equitable 
$(q, 1)$-tree-coloring for every $q\geq p.$ 
In this paper, we obtain a polynomial time criterion 
to decide if a complete bipartite graph has 
an equitable $(q,t)$-tree-coloring. 
Nevertheless, deciding if a graph $G$ in general 
has  an equitable $(q,t)$-tree-coloring is NP-complete.

\section{$(q,t)$-tree-coloring on graphs} 
Assume $G$ is a subgraph of $H,$ and $f$ is a coloring of $H,$ 
we let $f_G$ denote a coloring of $G$ defined by 
$f_G(v) = f(v)$ for each vertex $v$  in $G.$  

\begin{lem}\label{NPt}
Assume $t$ in a positive integer. 
Let $H$ be obtained from an $n$-vertex graph $G$ 
and $nt$ copies of $K_{2q-1}$ 
by joining each vertex of $G$ to $t$ copies of $K_{2q-1}$ 
where each copy of $K_{2q-1}$ is joining to exactly one 
vertex of $G.$ 
A graph $G$ has a proper $q$-coloring if and only if 
$H$ has a $(q,t)$-tree-coloring. 
\end{lem} 
\begin{pr} 
Observe that each color class resulted from $(q,\infty)$-tree-coloring 
of $K_{2q}$ has size 2. Thus such a coloring is unique 
(up to isomorphism).

Necessity. 
Color  $G$ by a proper $q$-coloring. 
We can extend a $(q,t)$-tree-coloring to each copy of $K_{2q-1}$ 
joining with a vertex in $G$ by the above observation. 
Moreover, there is exactly one vertex in each copy 
of $K_{2q-1}$ that has the same color with a vertex in $G.$ 
Thus this coloring is a $(q,t)$-tree-coloring of $H.$ 

Sufficiency. Assume $H$ has a $(q,t)$-tree-coloring $f.$  
By the above observation, each vertex $v$ in $G$ has exactly one adjacent vertex 
with the same color in each of  $t$ corresponding copies of $K_{2q-1}.$ 
Consequently, $v$ has no neighbors with the same color in $G.$ 
Thus $f_G$ is  a proper $q$-coloring. 
\end{pr}

\begin{lem}\label{NPi}
A graph $G$ has a proper $q$-coloring if and only if 
$H := G \vee K_q$ has a $(q,\infty)$-tree-coloring. 
\end{lem} 

\begin{pr}
Necessity. 
Color  $G$ by a proper $q$-coloring.
We can extend a $(q,\infty)$-tree-coloring to $H$ by coloring 
all vertices of $K_q$ by distinct colors. 
One can check that each color class induces $K_{1,m}$ for some $m.$ 
Thus we obtain a desired coloring. 

Sufficiency. 
For a coloring $f$ of $H,$ define $S(f)$ to be 
a number of color classes $V_i$s  
in which $V_i \cap V(G)$ is an independent set. 
Choose $f$ among $(q,\infty)$-tree-colorings of $H$ 
with the largest $S(f).$ 
Note that if $S(f) = q,$ then $f_G$ is a proper $q$-coloring of $G.$ 
Suppose to the contrary that $S(f) < q.$ 
Then there is a color class $V_1$ such that $V_1 \cap V(G)$ is not independent. 
Note that $V_1$ cannot contain a vertex from $K_q,$ otherwise $V_1$ induces a cycle subgraph. 
Since we color this $K_q$ by $q-1$ remaining colors, there is a color class $V_2 = \{a, b\} \subseteq V(K_q).$  
 
Since $V_1$ induces a forest, we can partition $V_1$ into 2 independent sets $X$ and $Y.$ 
Let a coloring $g$ be obtained from $f$ by changing $V_1$ into $X \cup \{a\}$ 
and $V_2$ into $Y \cup \{b\}$ while other color classes remain the same. 
Thus we obtain a  $(q,\infty)$-tree-coloring $g$ of $H$ with $S(g) = S(f)+1.$ 
This contradiction completes the proof.
\end{pr}

\begin{cor}
Let $q \geq \max \{3,t\}.$  
The problems of determining if a graph $G$ has a $(q,t)$-tree-coloring, 
a $(q,\infty)$-tree-coloring, 
an equitable  $(q,t)$-tree-coloring, 
or an equitable $(q,\infty)$-tree-coloring are NP-complete. 
\end{cor}

\begin{pr}
It is known \cite{GareyJohnson} that determining if a planar graph $G$ with
maximum vertex degree 4 is 3-colorable is NP-complete. 
Thus deciding if $G \vee K_{q-3}$ 
is $q$-colorable where $q \geq 4$ is also NP-complete. 
Using Lemmas \ref{NPt} and \ref{NPi}, we have that problems of 
determining if a graph has a $(q,t)$-tree-coloring or 
a $(q,\infty)$-tree-coloring for $q \geq \max\{3,t\}$ are NP-complete. 

Let $H$ be obtained from an $n$-vertex graph $G$ 
by adding $qn$ isolated vertices. 
Then $G$ has  a $(q,t)$-tree-coloring (respectively, 
a $(q,\infty)$-tree-coloring) if and only if $H$ 
has  an equitable $(q,t)$-tree-coloring (respectively, 
an equitable $(q,\infty)$-tree-coloring). 
This completes the proof. 
\end{pr}

\section{$(q,t)$-tree-coloring on bipartite graphs} 

Let $V_1,\ldots, V_q$ be color classes from   a $q$-coloring $c$ 
(not needed to be proper) of $K_{m,n},$  
$a =\lfloor (m+n)/q \rfloor,$ and let 
\begin{displaymath}
 \begin{array} {lr}
 c(X_1)=\{V_i | \; |V_i \cap X|= a+1, |V_i \cap Y|= 0\}, &
 c(X_2)=\{V_i | \;|V_i \cap X|= a, |V_i \cap Y|= 0\},\\
 c(Y_1)=\{V_i | \;|V_i \cap Y|= a+1, |V_i \cap X|= 0\}, &
 c(Y_2)=\{V_i | \;|V_i \cap Y|= a, |V_i \cap X|= 0\},\\
 c(X'_1)=\{V_i | \;|V_i \cap X|= a, |V_i \cap Y|= 1\}, &
 c(X'_2)=\{V_i | \;|V_i \cap X|= a-1, |V_i \cap Y|= 1\},\\ 
 c(Y'_1)=\{V_i | \;|V_i \cap Y|= a, |V_i \cap X|= 1\}, &
 c(Y'_2)=\{V_i | \;|V_i \cap Y|= a-1, |V_i \cap X|= 1\}.
 \end{array}
 \end{displaymath} 
 We have the following lemma. 
 
\begin{lem} \label{eqp} 
Let $a =\lfloor (m+n)/q \rfloor$ and $m+n =qa+r.$    
There are nonnegative integers $k_1,\ldots,k_4$ satisfying 
\begin{eqnarray}\label{eq1}
k_1 +k_2+k_3+k_4&=q, \nonumber \\
k_1(a+1)+k_2a &=m, \\
 k_3(a+1)+k_4a &=n \nonumber
\end{eqnarray} 
if and only if $K_{m,n}$ has a proper equitable $q$-coloring 
such that $|c(X_1)|=k_1, |c(X_2)|=k_2,$  
$|c(Y_1)|=k_3,$ and $|c(Y_2)|=k_4.$ 

Additionally, we have $k_1+k_3=r$ and $k_2+k_4=q-r$ for the above values. 
\end{lem}  
 
\begin{pr} 
If $c$ is a proper equitable $q$-coloring,  
then each color class is contained in $X$ or $Y$ 
and each has size $a$ or $a+1.$  
Thus $|c(X_1)|+|c(X_2)| +|c(Y_1)|+|c(Y_2)|=q,$ 
$|c(X_1)|(a+1)+ |c(X_2)|a =m,$ 
and  $|c(Y_1)|(a+1)+ |c(Y_2)|a =n.$  

Conversely, assume nonnegative integers $k_1,k_2,k_3,k_4$ 
satisfy all of these equations. 
It is straightforward to construct a proper equitable $q$-coloring $c$ 
with $|c(X_1)|=k_1, |c(X_2)|=k_2,$  
$|c(Y_1)|=k_3,$ and $|c(Y_2)|=k_4.$ Thus the converse holds. 

Consider $k_1,k_2,k_3,$ and $k_4.$ 
From $qa+r=m+n= (k_1+k_3)(a+1)+(k_2+k_4)a= (k_1+k_2+k_3+k_4)a+k_1+k_3=qa+(k_1+k_3),$ 
we have  $k_1+k_3=r$ and $k_2+k_4=q-r.$
\end{pr}
 
Let $a =\lfloor (m+n)/q \rfloor.$ 
Each color class $V_i$ from  
an equitable $(q,a)$-tree-coloring of $K_{m,n}$  
has size  $a$ or $a+1.$  
Moreover,  $|V_i \cap X| \leq 1$ or $|V_i \cap Y| \leq 1,$  
otherwise $V_i$ induces a subgraph $C_4.$  
Using these facts, Wu, Zhang, and Li obtained the analogous result 
for an equitable $(q,a)$-tree-coloring.

\begin{lem} \label{eqt} \cite{WZL13}
Let $a =\lfloor (m+n)/q \rfloor$ and $t \geq a.$   
 There are nonnegative integers $k_1,\ldots,k_8$ satisfying 
\begin{eqnarray}\label{eq2}
 k_1 +\cdots+k_8&=q, \nonumber\\ 
 k_1(a+1)+k_2a +k_5a+ k_6(a-1)+k_7+k_8 &=m, \\
 k_3(a+1)+k_4a +k_7a+ k_8(a-1)+k_5+k_6&=n \nonumber  
\end{eqnarray}
if and only if $K_{m,n}$ has an equitable $(q,t)$-tree-coloring $c$ 
such that $|c(X_1)|=k_1, |c(X_2)|=k_2,$  
$|c(Y_1)|=k_3, |c(Y_2)|=k_4,$ 
$|c(X'_1)|=k_5, |c(X'_2)|=k_6,$ 
$|c(Y'_1)|=k_7,$ and $|c(Y'_2)| =k_8.$ 

Additionally, we have $k_1+k_3+k_5+k_7=r$ and $k_2+k_4+k_6+k_8=q-r$ for the above values. 
\end{lem} 





\noindent {\bf Condition A.} Let $m+n=qa+r$ where $0\leq r \leq q-1.$ 
We call $(m,n)$ satisfies Condition (A) 
if one of the following holds; \\
(i) $r=0$ and $m$ is divisible by $a,$\\ 
(ii) $r \geq 1, r(a+1) \geq m,$ and 
$\min \{\lceil m/(a+1) \rceil, q-r \} \geq 
(a+1) \lceil m/(a+1) \rceil  -m,$ or \\
(iii) $r \geq 1, r(a+1) < m,$ and 
$\min \{r, q -\lceil (m-r)/a \rceil \} \geq 
a \lceil (m-r)/a \rceil  +r-m.$

\begin{thm}\label{thmA} 
$K_{m,n}$ has a proper equitable $q$-coloring if and only if 
$(m,n)$ or $(n,m)$ satisfies Condition A. 
\end{thm}
\begin{pr} 
NECESSITY. Assume $c$ is a proper equitable $q$-coloring 
of $K_{m,n}.$  
WLOG, assume $|c(X_1)| \geq |c(Y_1)|.$ 
We prove that $(m,n)$ satisfies Condition A. 

CASE 1. $r=0.$ 
Then every color class has size $a.$ 
Consequently, $ |c(X_2)|a=m.$  
Thus $(m,n)$ satisfies Condition A(i). 

CASE 2. $r\geq 1.$ 
Choose such a coloring $c$ with largest $|c(X_1)|.$ 
Suppose to the contrary that $|c(X_2)| \geq a+1$ 
and $|c(Y_1)| \geq a.$ 
Define $k_1= |c(X_1)| + a, k_2=|c(X_2)|-a-1, k_3=|c(Y_1)|-a,
k_4= |c(Y_2)|+a+1.$ 
By Theorem \ref{eqp}, there is a proper equitable $q$-coloring $c'$ of 
$K_{m,n}$ with $|c'(X_1)| = |c(X_1)|+a$ 
which contradicts the choice of $c.$ 
Thus $|c(X_2)| \leq a$ or $|c(Y_1)| \leq a-1.$ 

SUBCASE 2.1. $|c(X_2)| \leq |c(Y_1)|.$ 
Since $|c(X_1)|(a+1)+|c(X_2)|a=m,$ we have 
$(|(c(X_1)|+|c(X_2)|-1)(a+1)\leq  m \leq (|c(X_1)|+|c(X_2)|)(a+1).$ 

If $(|c(X_1)|+|c(X_2)|-1)(a+1) <  m,$ 
then $\lceil m/(a+1) \rceil =  |c(X_1)|+|c(X_2)|.$ 
Consequently, 
$\min \{\lceil m/(a+1) \rceil, q-r \}  
= \min\{|c(X_1)|+|c(X_2)|,  |c(X_2)|+|c(Y_2)|\}
\geq |c(X_2)|=(a+1) \lceil m/(a+1) \rceil  -m.$ 

If $(|c(X_1)|+|c(X_2)|-1)(a+1)=  m,$ 
then   $\min \{\lceil m/(a+1) \rceil, q-r \}  
= \min\{|c(X_1)|+|c(X_2)|-1,  |c(X_2)|+|c(Y_2)|\}
\geq 0=(a+1) \lceil m/(a+1) \rceil  -m.$  
Thus $(m,n)$ satisfies Condition A(ii). 

SUBCASE 2.2. $|c(X_2)| > |c(Y_1)|.$ 
Then $|c(Y_1)| \leq a-1.$ 
Combining with 
$m-r = |c(X_1)|(a+1)+|c(X_2)|a - (|c(X_1)|+|c(Y_1)|)
=(|c(X_1)|+|c(X_2)|)a-|c(Y_1)|,$ 
we have $\lceil (m-r)/a \rceil = |c(X_1)|+|c(X_2)|.$ 
Consequently,  $q -\lceil (m-r)/a \rceil 
=|c(X_1)|+|c(X_2)|+|c(Y_1)|+|c(Y_2)| -(|c(X_1)|+|c(X_2)|) =|c(Y_1)|+|c(Y_2)|.$ 
Moreover, $a \lceil (m-r)/a \rceil  -m+r 
= (|c(X_1)|+|c(X_2)|)a - |c(X_1)|(a+1)-|c(X_2)|a  + |c(X_1)|+|c(Y_1)|
=|c(Y_1)|.$  
 
Consequently, $\min \{r, q -\lceil (m-r)/a \rceil \} 
= \{ |c(X_1)|+|c(Y_1)|, |c(Y_1)|+|c(Y_2)|\}
\geq |c(Y_1)| = a \lceil (m-r)/a \rceil  -m+r.$ 
Thus $(m,n)$ satisfies Condition A(iii).

SUFFICIENCY. 

CASE 1.  $r=0$ and $m$ is divisible by $a.$ 
We can construct a proper equitable coloring $c$ 
with $|c(X_1)|=m/a, |c(Y_1)|=n/a,$ and $|c(X_2)|=|c(Y_2)|=0.$ 

CASE 2. $r \geq 1, r(a+1) \geq m,$ and 
$\min \{\lceil m/(a+1) \rceil, q-r \} \geq 
(a+1) \lceil m/(a+1) \rceil  -m.$ 

Choose $k_1=m-a\lceil m/(a+1) \rceil,$ 
$k_2= (a+1)\lceil m/(a+1) \rceil-m,$
$k_3=r+a\lceil m/(a+1) \rceil -m,$ and 
$k_4=m+q-r-(a+1)\lceil m/(a+1) \rceil.$ 

By assumption,  $k_1 = \lceil m/(a+1) \rceil 
- ((a+1)\lceil m/(a+1) \rceil-m) \geq 0,$ 
$k_3 \geq m/(a+1)+a\lceil m/(a+1) \rceil -m \geq 0,$  
and 
$k_4= (q-r)-((a+1)\lceil m/(a+1) \rceil-m) \geq 0.$ 
Obviously, $k_2$ is also nonnegative. 

It is straightforward to check that 
equation system \ref{eq1} in Lemma \ref{eqp} is satisfied.
Thus $K_{m,n}$ has a proper equitable $q$-coloring.

CASE 3. $r \geq 1, r(a+1) < m,$ and 
$\min \{r, q -\lceil (m-r)/a \rceil \} \geq 
a \lceil (m-r)/a \rceil  +r-m.$ 

Choose $k_1= m-a\lceil (m-r)/a \rceil,$ 
$k_2= (a+1)\lceil (m-r)/a \rceil -m,$ 
$k_3=r+ a\lceil (m-r)/a \rceil -m,$ and 
$k_4= m+q-r-(a+1)\lceil (m-r)/a \rceil.$ 

By assumption, $k_1=r-(a \lceil (m-r)/a \rceil +r -m) \geq 0,$  
$k_2 \geq (a+1)(m-r)/a -m= (m-r(a+1))/a \geq 0$   
and 
$k_4= (q -\lceil (m-r)/a \rceil)  -(a \lceil (m-r)/a \rceil +r-m) \geq 0.$ 
Obviously, $k_3$ is also nonnegative. 

It is straightforward to check that 
equation system \ref{eq1} in Lemma \ref{eqp} is satisfied.
Thus $K_{m,n}$ has a proper equitable $q$-coloring. 
\end{pr}

\noindent {\bf Condition B.} Let $m+n=qa+r$ where $0\leq r \leq q-1.$ 
We call $(m,n)$ satisfies Condition B if one of the following holds; \\
(i) $r \geq 1, r(a+1) \geq m,$ and 
$q+ a\lfloor m/(a+1) \rfloor \geq m,$ \\
(ii) $r=0$ or $r \geq 1, r(a+1) < m,$ and 
$ q+r+ (a-1)\lfloor (m-r)/a \rfloor 
\geq m.$ 

\begin{thm}\label{thmB}  
$K_{m,n}$ has an equitable $(q,t)$-tree-coloring 
if and only if 
$(m,n)$ or $(n,m)$ satisfies Condition A or B. 
\end{thm}

\begin{pr} 
SUFFICIENCY. 
Assume that $(m,n)$ satisfies Condition B. 

CASE 1. $r \geq 1, r(a+1) \geq m$ and 
$q+ a\lfloor m/(a+1) \rfloor \geq m.$ 

Choose $k_1= \lfloor m/(a+1) \rfloor,$ 
$k_2=k_5=k_6=0,$ 
$k_3= \max \{0, r+a\lfloor m/(a+1) \rfloor-m\}, $
$k_4=q+a\lfloor m/(a+1) \rfloor -m -k_3,$ 
$k_7= r-\lfloor m/(a+1) \rfloor -k_3,$ and 
$k_8=m+k_3-r-a\lfloor m/(a+1) \rfloor.$ 

By assumption, $k_4$ and $k_7$ are nonnegative. 
Obviously, each remaining $k_i$ is nonnegative. 
It is straightforward to check that 
equation system \ref{eq2} in Lemma \ref{eqt} is satisfied.
Thus $K_{m,n}$ has an equitable $(q,t)$-tree-coloring. 

CASE 2. $r=0$ and 
$ q+ (a-1)\lfloor (m-r)/a \rfloor \geq m.$ 

Choose $k_2= \lfloor m/a \rfloor,$ 
$k_4= q+(a-1)\lfloor m/a \rfloor-m,$ 
$k_8= m- a\lfloor m/a \rfloor,$ and 
$k_1=k_3=k_5=k_6=k_7=0.$ 
By assumption,  $k_4$ is nonnegative. 
Obviously, each remaining $k_i$ is nonnegative.  
It is straightforward to check that 
equation system \ref{eq2} in Lemma \ref{eqt} is satisfied.
Thus $K_{m,n}$ has an equitable $(q,t)$-tree-coloring.

CASE 3. $r \geq 1, r(a+1) < m$ and 
$ q+r+ (a-1)\lfloor (m-r)/a \rfloor \geq m.$ 

Choose $k_1= r,$ 
$k_2 = \lfloor (m-r)/a \rfloor -r,$ 
$k_3=k_5=k_6=k_7=0,$ 
$k_4= 
 q+r+ (a-1)\lfloor (m-r)/a \rfloor - m$ 
and 
$k_8= m-r - a \lfloor (m-r)/a \rfloor.$

By assumption, $k_2$ and $k_4$ are nonnegative. 
Obviously, each remaining $k_i$ is nonnegative. 
It is straightforward to check that 
equation system \ref{eq2} in Lemma \ref{eqt} is satisfied.
Thus $K_{m,n}$ has an equitable $(q,t)$-tree-coloring. 

Combining with Theorem \ref{thmA}, we complete the proof.

NECESSITY. Suppose $K_{m,n}$ has 
an equitable $(q,t)$-tree-coloring $c.$ 
We prove that $(m,n)$ or $(n,m)$ satisfies Condition A or B. 
If $c$ is an equitable $q$-coloring, 
then $(m,n)$ or $(n,m)$ satisfies Condition A. 
Assume equitable coloring does not exist for $K_{m,n}.$  
Consider equitable $(q,t)$-tree-colorings $c$ of $K_{m,n}.$ 
Define $k_1=|c(X_1)|, k_2= |c(X_2)|,$  
$k_3=|c(Y_1)|, k_4=|c(Y_2)|,$ 
$k_5=|c(X'_1)|, k_6=|c(X'_2)|,$ 
$k_7=|c(Y'_1)|,$ and $k_8=|c(Y'_2)|.$ 
We now restrict to $c$ with the least $k_5+k_6+k_7+k_8.$ 

We claim $|c(X'_1)|+ |c(X'_2)|=0$ or $|c(Y'_1)|+ |c(Y'_2)|=0.$  
Suppose to the contrary,  let $b= \min \{ 1,|c(X'_1)|\},$ 
$d= \min \{ 1,|c(Y'_1)|\},$ and redefine 
$k_5=|c(X'_1)|-b, k_6=|c(X'_2)|+b-1,$ 
$k_7=|c(Y'_1)|-d,$ and $k_8=|c(Y'_2)|+d-1.$ 
By Theorem \ref{eqt}, we have an equitable $(q,t)$-tree-colorings 
of $K_{m,n}$ with $|c(X'_1)|+ |c(X'_2)|+|c(Y'_1)|+ |c(Y'_2)|-2$  
defective color classes. This contradicts the choice of $c.$  
Thus the claim holds. 
WLOG, we assume $|c(X'_1)|+ |c(X'_2)|=0$ and choose such $c$ 
with largest $|c(X_1)|.$  
Since $c$ is not a proper coloring by assumption, 
we have $|c(Y'_1)|+ |c(Y'_2)|\geq 1.$ 

CASE 1. $r \geq 1, r(a+1) \geq m.$ 
If $|c(X_1)|(a+1)=m,$ then we can obtain an equitable coloring easily. 
Thus $|c(X_1)|(a+1)<m.$  
Since $|c(X_1)|(a+1)+ 
(|c(X'_1)|+|c(Y_1)|+|c'(Y_1)|)(a+1)=r(a+1)\geq m$ and 
$|c(X'_1)|=0,$ we have $|c(Y_1)|+|c'(Y_1)|\geq 1$ 
by condition of the case. 

We claim $|c(X_2)|=0.$  
Suppose to the contrary. 
Redefine $k_1=|c(X_1)|+1$ and $k_2= |c(X_2)|-1.$ 
If $|c(Y'_1)| \geq 1,$ then redefine 
$k_4=|c(Y_2)|+1,$ and $k_7=|c(Y'_1)|-1,$ 
otherwise redefine 
$k_3=|c(Y_1)|-1, k_4=|c(Y_2)|+2,$ 
and $k_8=|c(Y'_2)|-1.$ 
We reach a contradiction by a smaller value of $k_5+k_6+k_7+k_8.$ 

We claim further $|c(X_1)|=\lfloor m/(a+1) \rfloor$ or equivalently 
$m- |c(X_1)|(a+1) = |c'(Y_1)|+|c'(Y_2)| \leq a.$ 
Suppose to the contrary. 
If $|c'(Y_1)|\geq 1,$ then let $h= \min \{ a,|c(Y'_2)|\}$  
and redefine 
$k_1=|c(X_1)|+1,$ 
$k_3=|c(Y_1)|+a-h, k_4=|c(Y_2)|+h,$  
$k_7=|c(Y'_1)|-a+h-1,$ and $k_8=|c(Y'_2)|-h,$ 
otherwise redefine  
$k_1=|c(X_1)|+1,$  
$k_3=|c(Y_1)|-1, k_4=|c(Y_2)|+a+1,$ 
and $k_8=|c(Y'_2)|-a-1.$ 
Again, we reach a contradiction by a smaller value of $k_5+k_6+k_7+k_8.$ 

From above, we have 
$m = |c(X_1)|(a+1)+|c(Y'_1)|+|c(Y'_2)|$ 
$\leq  |c(X_1)|(a+1)+q - |c(X_1)|$ 
$ = a|c(X_1)|+q = a\lfloor (m/(a+1)) \rfloor +q.$ 
Thus $(m,n)$ satisfies Condition B.

CASE 2. $r =0.$ 
If $|c(X_2)|a=m,$ then we can obtain an equitable coloring easily. 
Thus we assume $m$ is not divisible by $a.$ 

We claim further $|c(X_2)|=\lfloor m/a \rfloor$ or equivalently 
$m- |c(X_2)|a = |c'(Y_2)| \leq a-1.$ 
Suppose to the contrary. 
Redefine 
$k_2=|c(X_2)|+1,$ 
$k_4=|c(Y_2)|+a-1,$  
and $k_8=|c(Y'_2)|-a.$ 
We reach a contradiction by a smaller value of $k_5+k_6+k_7+k_8.$ 

From above, we have 
$q-\lfloor m/a \rfloor = q - |c(X_2)|= |c(Y_2)|+|c(Y'_2)|$ 
$\geq |c(Y'_2)| m-\lfloor m/a \rfloor a. $ 
That is $q+ (a-1)\lfloor m/a \rfloor \geq m.$ 
Thus $(m,n)$ satisfies Condition B.

CASE 3. $r \geq 1$ and  $r(a+1) < m.$

We claim that $m-|c(X_1)|(a+1)-|c(X_2)|a \leq a.$ 
Suppose to the contrary. 
If $|c'(Y_1)|\geq 1,$ then let $h= \min \{a+1,|c(Y'_1)|\}$  
and redefine 
$k_1=|c(X_1)|+1,$ 
$k_3=|c(Y_1)|+h-1, k_4=|c(Y_2)|+a-h+1,$  
$k_7=|c(Y'_1)|-h,$ and $k_8=|c(Y'_2)|+h-a-1,$ 
otherwise redefine  
$k_2=|c(X_2)|+1,$  
$k_4=|c(Y_2)|+a-1,$ 
and $k_8=|c(Y'_2)|-a.$ 
We reach a contradiction by a smaller value of $k_5+k_6+k_7+k_8.$  

We claim further that $|c(X_1)|=r.$ 
Combining with the previous claim, 
we have $|c(X_2)|=\lfloor (m-r(a+1))/a \rfloor a.$  
Suppose to the contrary.  
Then $|c(Y_1)|+|c(Y'_1)| \geq 1.$ 
From $m-|c(X_1)|(a+1)-|c(X_2)|a \leq a$ and condition of the case, 
we have $|c(X_2)|\geq 1.$ 
Since $c$ is not a proper coloring by assumption, 
$|c(Y'_1)|+|c(Y'_2)|\geq 1.$ 
Redefine $k_1=|c(X_1)|+1$ and 
$k_2=|c(X_2)|-1.$ 
If $|c'(Y_1)|\geq 1,$ then redefine 
$k_4=|c(Y_2)|+1,$ and 
$k_7=|c(Y'_1)|-1,$ 
otherwise redefine  
$k_3=|c(Y_1)|-1,$ 
$k_4=|c(Y_2)|+2,$ and 
$k_8=|c(Y'_2)|-1.$ 

Thus $m-|c(X_1)|(a+1)-|c(X_2)|a  \leq$ 
$m-r(a+1)-\lfloor (m-r(a+1))/a \rfloor a = |c(Y'_1)|+|c(Y'_2)|$
$\leq |c(Y_1)|+|c(Y_2)|+|c(Y'_1)|+|c(Y'_2)|$ 
$= q - c(X_1) - c(X_2)= q-r-\lfloor (m-r(a+1))/a \rfloor.$ 
That is 
$ q+r+ (a-1)\lfloor (m-r)/a \rfloor \geq m.$ 
Thus $(m,n)$ satisfies Condition B. 
\end{pr}

\end{document}